\documentclass[10pt]{article}

\usepackage{amsmath,amssymb,amscd,amsthm,makeidx,txfonts,graphicx,url,bm}
\usepackage{a4wide,pstricks,xypic}

\numberwithin{equation}{section}

\let\oldmarginpar\marginpar
\renewcommand\marginpar[1]{\-\oldmarginpar[\raggedleft\footnotesize #1]
{\raggedright\footnotesize #1}}

\makeindex

\xyoption{all}
\input{xypic}

\newtheorem{theorem}{Theorem}[section]

\newtheorem{corollary}[theorem]{Corollary}
\newtheorem{conjecture}[theorem]{Conjecture}

\theoremstyle{remark}

\theoremstyle{definition}

\def\beq{\begin{eqnarray}}
\def\eeq{\end{eqnarray}}
\def\bes{\begin{eqnarray*}}
\def\ees{\end{eqnarray*}}

\def\muhat{{\bm \mu}}

\def\C{\mathbb{C}}
\def\M{{\mathcal{M}}}
\def\calQ{{\mathcal{Q}}}

\def\x{\mathbf{x}}
\def\y{\mathbf{y}}

\def\v{\mathbf{v}}

\def\P{\mathcal{P}}
\def\H{\mathbb{H}}

\def\F{\mathbb{F}}
\def\Q{\mathbb{Q}}

\def\calC{{\mathcal C}}
\def\calO{{\mathcal O}}

\def\Z{\mathbb{Z}}

\def\K{\mathbb{K}}

\def\gl{{\mathfrak g\mathfrak l}}

\newcommand{\nc}{\newcommand}

\nc{\op}[1]{\mathop{\mathchoice{\mbox{\rm #1}}{\mbox{\rm #1}}
{\mbox{\rm \scriptsize #1}}{\mbox{\rm \tiny #1}}}\nolimits}
\nc{\al}{\alpha}

\nc{\ep}{\varepsilon} \nc{\ga}{\gamma} \nc{\Ga}{\Gamma}
\nc{\la}{\lambda} \nc{\La}{\Lambda} \nc{\si}{\sigma}
\nc{\Sig}{{\Gamma}} \nc{\Om}{\Omega} \nc{\om}{\omega}

\nc{\SL}{{\rm SL}} \nc{\GL}{{\rm GL}} \nc{\PGL}{{\rm PGL}}
\nc{\G}{{\rm G}}

\nc{\cpt}{{\op{cpt}}} \nc{\Dol}{{\op{Dol}}} \nc{\DR}{{\op{DR}}}
\nc{\B}{{\op{B}}} \nc{\Triv}{\op{Triv}} \nc{\Hod}{{\op{Hod}}}
\nc{\Log}{{\op{Log}}} \nc{\Exp}{{\op{Exp}}} \nc{\Est}{E_{\op{st}}}
\nc{\Hst}{H_{\op{st}}} \nc{\Left}[1]{\hbox{$\left#1\vbox to
  10.5pt{}\right.\nulldelimiterspace=0pt \mathsurround=0pt$}}
\nc{\Right}[1]{\hbox{$\left.\vbox to
  10.5pt{}\right#1\nulldelimiterspace=0pt \mathsurround=0pt$}}
\nc{\LEFT}[1]{\hbox{$\left#1\vbox to
  15.5pt{}\right.\nulldelimiterspace=0pt \mathsurround=0pt$}}
\nc{\RIGHT}[1]{\hbox{$\left.\vbox to
  15.5pt{}\right#1\nulldelimiterspace=0pt \mathsurround=0pt$}}

\nc{\bee}{{\bf E}} \nc{\bphi}{{\bf \Phi}}

\begin{document}

\title{Topology of character varieties and representations of quivers}

\author{ Tam\'as Hausel
\\ {\it University of Oxford}
\\ {\it
  University of Texas at Austin}
\\{\tt hausel@maths.ox.ac.uk} \and Emmanuel Letellier \\ {\it Universit\'e de Caen} \\ {\tt letellier.emmanuel@math.unicaen.fr}\and Fernando Rodriguez-Villegas \\
{\it University of Texas at Austin} \\ {\tt
  villegas@math.utexas.edu}\\
}

\pagestyle{myheadings}

\maketitle

\begin{abstract}

 In \cite{hausel-letellier-villegas} we
presented a conjecture generalizing the Cauchy formula for Macdonald
polynomials. This conjecture encodes the mixed Hodge polynomials of
the character varieties of representations of the fundamental group
of a punctured Riemann surface of genus $g$. We proved several
results which support this conjecture. Here we announce new results
which are consequences of those in
\cite{hausel-letellier-villegas}.
\end{abstract}

\section{Review of the results of \cite{hausel-letellier-villegas}}

\subsection{Cauchy function} Fix integers $g\geq 0$ and $k>0$. Let
$\x_1=\{x_{1,1},x_{1,2},...\},...,\x_k=\{x_{k,1},x_{k,2},...\}$ be
$k$ sets of infinitely many independent variables and let $\Lambda$
be the ring of functions separately symmetric in each set of
variables. Let $\P$ be the set of partitions. For $\lambda\in\P$,
let $\tilde{H}_{\lambda}(\x_i;q,t)\in\Lambda\otimes\Q(q,t)$ be the
\emph{Macdonald symmetric function} defined in
\cite[I.11]{garsia-haiman}.

Define the \emph{$k$-point genus $g$ Cauchy function}
$$
\Omega(z,w)=
\sum_{\lambda\in\P}
\mathcal{H}_{\lambda}(z,w)
\prod_{i=1}^k\tilde{H}_{\lambda}\big(\x_i;z^2,w^2\big).
$$
where
$$
\mathcal{H}_{\lambda}(z,w)
:=\prod\frac{(z^{2a+1}-w^{2l+1})^{2g}}{(z^{2a+2}-w^{2l})(z^{2a}-w^{2l+2})}
$$
is a $(z,w)$-deformation of the $(2g-2)$-th power of the standard
hook polynomial. Let $\Exp$ be the plethystic exponential and let
$\Log$ be its inverse~\cite[2.3]{hausel-letellier-villegas}. For
$\muhat=(\mu^1,...,\mu^k)\in\P^k$, let
$$\H_{\muhat}(z,w):=(z^2-1)(1-w^2)\left\langle\Log(\Omega(z,w),
h_{\muhat}\right\rangle$$where $h_{\mu}:=h_{\mu^1}(\x_1)\cdots
h_{\mu^k}(\x_k)\in\Lambda$ is the product of the complete symmetric
functions and $\langle,\rangle$ is the extended Hall pairing.

\subsection{Character and quiver varieties} Let $\M_{\muhat}$ and $\calQ_{\muhat}$ be
\emph{generic} character and quiver varieties corresponding to
$\muhat$ \cite[2.1,2.2]{hausel-letellier-villegas}. Namely, we let
$(\calC_1,...,\calC_k)$  be a \emph{generic} $k$-tuple of semisimple conjugacy
classes of $\GL_n(\C)$ of type $\muhat$, i.e., the coordinate
$\mu^i$ of $\muhat$ gives the multiplicities of the eigenvalues of
$C_i$. Then $\M_{\muhat}$ is the affine GIT quotient

\bes \M_{\muhat}:= \{ A_1,B_1,\dots,A_g,B_g \in \GL_n(\C), X_1\in
\calC_1,\dots, X_k\in \calC_k |\, (A_1,B_1) \cdots (A_g,B_g)
X_1\cdots X_k=I_n\}/\!/ \GL_n(\C), \ees where for two matrices
$A,B$, we denote by $(A,B)$ the commutator $ABA^{-1}B^{-1}$. Let
$(\calO_1,...,\calO_k)$ be a generic $k$-tuple of semisimple adjoint orbits of
$\gl_n(\C)$ of type $\muhat$. The quiver variety $\calQ_{\muhat}$ is
defined as the affine GIT quotient\vspace{.2cm}

\noindent $\calQ_{\muhat}:= \{ A_1,B_1,\dots,A_g,B_g \in \gl_n(\C),
C_1\in \calO_1,\dots, C_k\in \calO_k |\, [A_1,B_1]+ \cdots
+[A_g,B_g]+C_1+\cdots+C_k=O\}/\!/ \GL_n(\C).$\vspace{.2cm}

In \cite{hausel-letellier-villegas}, we proved that $\M_{\muhat}$
and $\calQ_{\muhat}$ are non-singular algebraic varieties  of pure
dimension $d_{\muhat}=n^2(2g-2+k)-\sum_{i,j}(\mu^i_j)^2+2$.

Let $H_c(\M_{\muhat};x,y,t):=\sum_{i,j,k}
h_c^{i,j;k}(\M_{\muhat})x^iy^jt^k$ be the compactly supported mixed
Hodge polynomial. It is a common deformation of the compactly
supported Poincar\'e polynomial
$P_c(\M_{\muhat};t)=H_c(\M_{\muhat};1,1,t)$ and the so-called
$E$-polynomial $E(\M_{\muhat};x,y)=H_c(\M_{\muhat};x,y,-1)$. We have
the following conjecture~\cite[Conjecture
1.1.1]{hausel-letellier-villegas}:
\begin{conjecture}
\label{main-conjecture}
 The polynomial $H_c(\M_{\muhat};x,y,t)$ depends only
on $xy$ and $t$. If we let
$H_c(\M_{\muhat};q,t)=H_c(\M_{\muhat};\sqrt{q},\sqrt{q},t)$ then
\begin{equation}
H_c(\M_{\muhat};q,t)=
(t\sqrt{q})^{d_{\muhat}}\H_{\muhat}
\left(-\frac{1}{\sqrt{q}},t\sqrt{q}\right).\label{conj}
\end{equation}
\end{conjecture}
This conjecture implies the following one:
\begin{conjecture}[Curious Poincar\'e duality]
$$
H_c\left(\M_{\muhat};\frac{1}{qt^2},t\right)
=(qt)^{-d_{\muhat}}H_c(\M_{\muhat};q,t).
$$
\label{curiousduality}
\end{conjecture}

The two following theorems are proved in~\cite{hausel-letellier-villegas}:

\begin{theorem} The $E$-polynomial $E(\M_{\muhat};x,y)$ depends only
on $xy$ and if we let
$E(\M_{\muhat};q)=E(\M_{\muhat};\sqrt{q},\sqrt{q})$, we have
\begin{equation}E(\M_{\muhat};q)=q^{\frac{1}{2}d_{\muhat}}\H_{\muhat}\left(\frac{1}{\sqrt{q}},\sqrt{q}\right).\label{theo1.1}\end{equation}
\label{epolychar}\end{theorem}

As a corollary we get a consequence of the curious Poincar\'e
duality Conjecture \ref{curiousduality}:

\begin{corollary} The $E$-polynomial is \emph{palindromic}.

$$E(\mathcal{M}_{\muhat};q)=q^{d_{\muhat}}E(\mathcal{M}_{\muhat};q^{-1})=\sum_i\left(\sum_k(-1)^kh^{i,i;k}(\M_{\muhat})\right)q^i.$$\label{curious}

\end{corollary}

We say that $\muhat$ is \emph{indivisible} if the gcd of all the
parts of the partitions $\mu^1,...,\mu^k$ is equal to $1$. It is
possible to choose $k$ generic semisimple adjoint orbits of type
$\muhat$ if and only if $\muhat$ is indivisible \cite[Lemma
2.2.2]{hausel-letellier-villegas}.

\begin{theorem} For $\muhat$ indivisible, the mixed Hodge structure
on $H_c^*(\calQ_{\muhat},\C)$ is pure. If we let
$E(\calQ_{\muhat};q)=E(\calQ_{\muhat};\sqrt{q},\sqrt{q})$, then
\begin{equation}P_c(\calQ_{\muhat};\sqrt{q})=
E(\calQ_{\muhat};q)=q^{\frac{1}{2}d_{\muhat}}\H_{\muhat}(0,\sqrt{q}).
\label{theo1.2}
\end{equation}
\end{theorem}

Note that Formula \eqref{theo1.1} is the specialization $t\mapsto -1$
of Formula \eqref{conj}. Assuming Conjecture \ref{main-conjecture},
Formula \eqref{theo1.2} implies that the $i$-th Betti number of
$\calQ_{\muhat}$ equals the dimension of the $i$-th piece of the pure part
of the cohomology of $\M_{\muhat}$, namely,
$\sum_ih_c^{i,i;2i}(\M_\muhat)$. Furthermore, when $g=0$, the first
author conjectures \cite{hausel2} that there is an isomorphism between
the pure part of $H_c^i(\M_{\muhat},\C)$ and
$H_c^i(\calQ_{\muhat},\C)$ induced by the Riemann-Hilbert monodromy
map $\calQ_{\muhat}\rightarrow \M_{\muhat}$. This would give a
geometric interpretation of Theorem \ref{theo1.2} in this case.

\subsection{Multiplicities in tensor products}

Given $\muhat=(\mu^1,...,\mu^k)\in\P^k$, we can choose a \emph{generic} $k$-tuple $(R_1,...,R_k)$ of
semisimple irreducible complex characters of
$\GL_n(\F_q)$ where $\F_q$ is a finite field with $q$ elements
\cite{hausel-letellier-villegas}. We also denote by
$\Lambda:\GL_n(\F_q)\rightarrow \C$ the character $h\mapsto
q^{g\,\text{dim}\,Z(h)}$ where $Z(h)$ is the centralizer of $h$ in $\GL_n(\overline{\F}_q)$. Then we have \cite[6.1.1]{hausel-letellier-villegas}:
\begin{theorem} $$\left\langle\Lambda\otimes
R_{\muhat},1\right\rangle=\H_{\muhat}(0,\sqrt{q})$$where
$R_{\muhat}=\bigotimes_{i=1}^kR_i$. \label{multi}\end{theorem}

\section{Absolutely indecomposable representations}

Let $\mathbf{s}=(s_1,...,s_k)\in\big(\Z_{\geq 0}\big)^k$. Let $\Gamma$ be
the comet-shaped quiver with $g$ loops on the central vertex
represented as below:

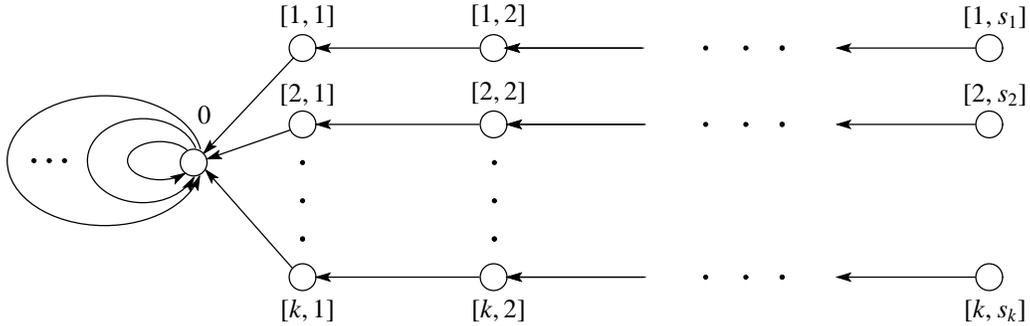
\begin{center}
\unitlength 0.1in
\begin{picture}( 52.1000, 15.4500)(  4.0000,-17.0000)
%
\special{pn 8}%
\special{ar 1376 1010 70 70  0.0000000 6.2831853}%
%
\special{pn 8}%
\special{ar 1946 410 70 70  0.0000000 6.2831853}%
%
\special{pn 8}%
\special{ar 2946 410 70 70  0.0000000 6.2831853}%
%
\special{pn 8}%
\special{ar 5540 410 70 70  0.0000000 6.2831853}%
%
\special{pn 8}%
\special{ar 1946 810 70 70  0.0000000 6.2831853}%
%
\special{pn 8}%
\special{ar 2946 810 70 70  0.0000000 6.2831853}%
%
\special{pn 8}%
\special{ar 5540 810 70 70  0.0000000 6.2831853}%
%
\special{pn 8}%
\special{ar 1946 1610 70 70  0.0000000 6.2831853}%
%
\special{pn 8}%
\special{ar 2946 1610 70 70  0.0000000 6.2831853}%
%
\special{pn 8}%
\special{ar 5540 1610 70 70  0.0000000 6.2831853}%
%
\special{pn 8}%
\special{pa 1890 1560}%
\special{pa 1440 1050}%
\special{fp}%
\special{sh 1}%
\special{pa 1440 1050}%
\special{pa 1470 1114}%
\special{pa 1476 1090}%
\special{pa 1500 1088}%
\special{pa 1440 1050}%
\special{fp}%
%
\special{pn 8}%
\special{pa 2870 410}%
\special{pa 2020 410}%
\special{fp}%
\special{sh 1}%
\special{pa 2020 410}%
\special{pa 2088 430}%
\special{pa 2074 410}%
\special{pa 2088 390}%
\special{pa 2020 410}%
\special{fp}%
%
\special{pn 8}%
\special{pa 3720 410}%
\special{pa 3010 410}%
\special{fp}%
\special{sh 1}%
\special{pa 3010 410}%
\special{pa 3078 430}%
\special{pa 3064 410}%
\special{pa 3078 390}%
\special{pa 3010 410}%
\special{fp}%
\special{pa 3730 410}%
\special{pa 3010 410}%
\special{fp}%
\special{sh 1}%
\special{pa 3010 410}%
\special{pa 3078 430}%
\special{pa 3064 410}%
\special{pa 3078 390}%
\special{pa 3010 410}%
\special{fp}%
%
\special{pn 8}%
\special{pa 2870 810}%
\special{pa 2020 810}%
\special{fp}%
\special{sh 1}%
\special{pa 2020 810}%
\special{pa 2088 830}%
\special{pa 2074 810}%
\special{pa 2088 790}%
\special{pa 2020 810}%
\special{fp}%
%
\special{pn 8}%
\special{pa 2870 1610}%
\special{pa 2020 1610}%
\special{fp}%
\special{sh 1}%
\special{pa 2020 1610}%
\special{pa 2088 1630}%
\special{pa 2074 1610}%
\special{pa 2088 1590}%
\special{pa 2020 1610}%
\special{fp}%
%
\special{pn 8}%
\special{pa 3730 810}%
\special{pa 3020 810}%
\special{fp}%
\special{sh 1}%
\special{pa 3020 810}%
\special{pa 3088 830}%
\special{pa 3074 810}%
\special{pa 3088 790}%
\special{pa 3020 810}%
\special{fp}%
\special{pa 3740 810}%
\special{pa 3020 810}%
\special{fp}%
\special{sh 1}%
\special{pa 3020 810}%
\special{pa 3088 830}%
\special{pa 3074 810}%
\special{pa 3088 790}%
\special{pa 3020 810}%
\special{fp}%
%
\special{pn 8}%
\special{pa 3730 1610}%
\special{pa 3020 1610}%
\special{fp}%
\special{sh 1}%
\special{pa 3020 1610}%
\special{pa 3088 1630}%
\special{pa 3074 1610}%
\special{pa 3088 1590}%
\special{pa 3020 1610}%
\special{fp}%
\special{pa 3740 1610}%
\special{pa 3020 1610}%
\special{fp}%
\special{sh 1}%
\special{pa 3020 1610}%
\special{pa 3088 1630}%
\special{pa 3074 1610}%
\special{pa 3088 1590}%
\special{pa 3020 1610}%
\special{fp}%
%
\special{pn 8}%
\special{pa 5466 410}%
\special{pa 4746 410}%
\special{fp}%
\special{sh 1}%
\special{pa 4746 410}%
\special{pa 4812 430}%
\special{pa 4798 410}%
\special{pa 4812 390}%
\special{pa 4746 410}%
\special{fp}%
%
\special{pn 8}%
\special{pa 5466 810}%
\special{pa 4746 810}%
\special{fp}%
\special{sh 1}%
\special{pa 4746 810}%
\special{pa 4812 830}%
\special{pa 4798 810}%
\special{pa 4812 790}%
\special{pa 4746 810}%
\special{fp}%
%
\special{pn 8}%
\special{pa 5466 1610}%
\special{pa 4746 1610}%
\special{fp}%
\special{sh 1}%
\special{pa 4746 1610}%
\special{pa 4812 1630}%
\special{pa 4798 1610}%
\special{pa 4812 1590}%
\special{pa 4746 1610}%
\special{fp}%
%
\special{pn 8}%
\special{pa 1880 840}%
\special{pa 1450 990}%
\special{fp}%
\special{sh 1}%
\special{pa 1450 990}%
\special{pa 1520 988}%
\special{pa 1500 972}%
\special{pa 1506 950}%
\special{pa 1450 990}%
\special{fp}%
%
\special{pn 8}%
\special{pa 1900 460}%
\special{pa 1430 960}%
\special{fp}%
\special{sh 1}%
\special{pa 1430 960}%
\special{pa 1490 926}%
\special{pa 1468 922}%
\special{pa 1462 898}%
\special{pa 1430 960}%
\special{fp}%
%
\special{pn 8}%
\special{sh 1}%
\special{ar 1946 1010 10 10 0  6.28318530717959E+0000}%
\special{sh 1}%
\special{ar 1946 1210 10 10 0  6.28318530717959E+0000}%
\special{sh 1}%
\special{ar 1946 1410 10 10 0  6.28318530717959E+0000}%
\special{sh 1}%
\special{ar 1946 1410 10 10 0  6.28318530717959E+0000}%
%
\special{pn 8}%
\special{sh 1}%
\special{ar 4056 410 10 10 0  6.28318530717959E+0000}%
\special{sh 1}%
\special{ar 4266 410 10 10 0  6.28318530717959E+0000}%
\special{sh 1}%
\special{ar 4456 410 10 10 0  6.28318530717959E+0000}%
\special{sh 1}%
\special{ar 4456 410 10 10 0  6.28318530717959E+0000}%
%
\special{pn 8}%
\special{sh 1}%
\special{ar 4056 810 10 10 0  6.28318530717959E+0000}%
\special{sh 1}%
\special{ar 4266 810 10 10 0  6.28318530717959E+0000}%
\special{sh 1}%
\special{ar 4456 810 10 10 0  6.28318530717959E+0000}%
\special{sh 1}%
\special{ar 4456 810 10 10 0  6.28318530717959E+0000}%
%
\special{pn 8}%
\special{sh 1}%
\special{ar 4056 1610 10 10 0  6.28318530717959E+0000}%
\special{sh 1}%
\special{ar 4266 1610 10 10 0  6.28318530717959E+0000}%
\special{sh 1}%
\special{ar 4456 1610 10 10 0  6.28318530717959E+0000}%
\special{sh 1}%
\special{ar 4456 1610 10 10 0  6.28318530717959E+0000}%
\put(19.7000,-2.4500){\makebox(0,0){$[1,1]$}}%
\put(29.7000,-2.4000){\makebox(0,0){$[1,2]$}}%
\put(55.7000,-2.5000){\makebox(0,0){$[1,s_1]$}}%
\put(19.7000,-6.5500){\makebox(0,0){$[2,1]$}}%
\put(29.7000,-6.4500){\makebox(0,0){$[2,2]$}}%
\put(55.7000,-6.5500){\makebox(0,0){$[2,s_2]$}}%
\put(19.7000,-17.8500){\makebox(0,0){$[k,1]$}}%
\put(29.7000,-17.8500){\makebox(0,0){$[k,2]$}}%
\put(55.7000,-17.8500){\makebox(0,0){$[k,s_k]$}}%
\put(14.3000,-7.6000){\makebox(0,0){$0$}}%
\special{pn 8}%
\special{sh 1}%
\special{ar 2950 1010 10 10 0  6.28318530717959E+0000}%
\special{sh 1}%
\special{ar 2950 1210 10 10 0  6.28318530717959E+0000}%
\special{sh 1}%
\special{ar 2950 1410 10 10 0  6.28318530717959E+0000}%
\special{sh 1}%
\special{ar 2950 1410 10 10 0  6.28318530717959E+0000}%
\special{pn 8}%
\special{ar 1110 1000 290 220  0.4187469 5.9693013}%
\special{pn 8}%
\special{pa 1368 1102}%
\special{pa 1376 1090}%
\special{fp}%
\special{sh 1}%
\special{pa 1376 1090}%
\special{pa 1324 1138}%
\special{pa 1348 1136}%
\special{pa 1360 1158}%
\special{pa 1376 1090}%
\special{fp}%
\special{pn 8}%
\special{ar 910 1000 510 340  0.2464396 6.0978374}%
\special{pn 8}%
\special{pa 1400 1096}%
\special{pa 1406 1084}%
\special{fp}%
\special{sh 1}%
\special{pa 1406 1084}%
\special{pa 1362 1138}%
\special{pa 1384 1132}%
\special{pa 1398 1152}%
\special{pa 1406 1084}%
\special{fp}%
\special{pn 8}%
\special{sh 1}%
\special{ar 540 1000 10 10 0  6.28318530717959E+0000}%
\special{sh 1}%
\special{ar 620 1000 10 10 0  6.28318530717959E+0000}%
\special{sh 1}%
\special{ar 700 1000 10 10 0  6.28318530717959E+0000}%
\special{pn 8}%
\special{ar 1200 1000 170 100  0.7298997 5.6860086}%
\special{pn 8}%
\special{pa 1314 1076}%
\special{pa 1328 1068}%
\special{fp}%
\special{sh 1}%
\special{pa 1328 1068}%
\special{pa 1260 1084}%
\special{pa 1282 1094}%
\special{pa 1280 1118}%
\special{pa 1328 1068}%
\special{fp}%
\end{picture}%
\end{center}\vspace{.2cm}

\noindent Let $I=\{0\}\cup\{[i,j]\}_{1\leq i\leq k,1\leq j\leq s_i}$
denote the set of vertices and let $\Omega$ be the set of arrows.
For $\gamma\in\Omega$, we denote by $h(\gamma)\in I$ the head of
$\gamma$ and $t(\gamma)\in I$ the tail of $\gamma$. A
\emph{dimension vector} for $\Gamma$ is a collection of non-negative
integers $\v=\{v_i\}_{i\in I}$ and a \emph{representation} $\varphi$
of $\Gamma$ of dimension $\v$ over a field $\K$ is a collection of
$\K$-vector spaces $\{V_i\}_{i\in I}$ with $\text{dim}\,V_i=v_i$
together with a collection of $\K$-linear maps
$\{\varphi_{\gamma}:V_{t(\gamma)}\rightarrow
V_{h(\gamma)}\}_{\gamma\in\Omega}$. We denote by ${\rm
Rep}_{\K}(\Gamma,\v)$ the $\K$-vector space of all representations
of $\Gamma$ over $\K$ of dimension vector $\v$. We also denote by
${\rm Rep}_{\K}^*(\Gamma,\v)$ the subset of representations
$\varphi\in{\rm Rep}_{\K}(\Gamma,\v)$ such that $\varphi_{\gamma}$
is injective for all $\gamma\in\Omega$ such that $t(\gamma)$ is not
the central vertex $0$.

Assume from now that $\K$ is a finite field $\F_q$. We denote by
${\rm Rep}^{a.i}_{\K}(\Gamma,\v)$ the set of absolutely
indecomposable representations in ${\rm Rep}_{\K}(\Gamma,\v)$. We
also assume that $v_0\neq 0$. Under this assumption, note that ${\rm
Rep}^{a.i}_{\K}(\Gamma,\v)\subset{\rm Rep}_{\K}^*(\Gamma,\v)$. We
may assume that $v_0\geq v_{[i,1]}\geq...\geq v_{i,s_i}$ for all
$i\in\{1,...,k\}$ since otherwise ${\rm Rep}^{a.i}_{\K}(\Gamma,\v)$
is empty. For each $i$, take the strictly decreasing subsequence
$v_0> n_{i_1}>...>n_{i_r}$ of $v_0\geq v_{[i,1]}\geq...\geq
v_{i,s_i}$ of maximal length. This defines a partition
$\mu^i:=\mu^i_1+...+\mu^i_{r+1}$ of $v_0$ as follows:
$\mu^i_1=v_0-n_{i_1}, 
\mu^i_2=n_{i_1}-n_{i_2},..., \mu^i_{r+1}=n_{i_r}$. The dimension
vector $\v$ defines thus a unique multipartition
$\muhat=(\mu^1,...,\mu^k)\in\P^k$.  The number $A_{\muhat}(q)$ of
isomorphism classes in ${\rm Rep}^{a.i}_{\K}(\Gamma,\v)$ depends
only on $\muhat$ and not on the choice of $\v$.

We have the following theorem \cite{hausel-letellier-villegas2}:
\begin{theorem}
For any $\muhat\in\P^k$
$$
A_{\muhat}(q)=\H_{\muhat}(0,\sqrt{q}).
$$
\label{abs}
\end{theorem}
We know by a theorem of V. Kac that $A_{\muhat}(q)\in\Z[q]$, see
\cite{kac}. It is also conjectured in \cite{kac} that the coefficients
of $A_{\muhat}(q)$ are non-negative. Assuming Conjecture
\ref{main-conjecture}, Theorem \ref{abs} gives a cohomological
interpretation of $A_{\muhat}(q)$; indeed, it implies that
$A_{\muhat}(q)$ is the Poincar\'e polynomial of the pure part of the
cohomology of $\M_{\muhat}$, thus implying Kac's conjecture for
comet-shaped quivers. In particular, combining Conjectures
\ref{main-conjecture} and \ref{curiousduality} and Theorem \ref{abs}
we obtain the conjectural equality of the middle Betti number of
$\M_\muhat$ and $A_\muhat(1)$. These remarks can be compared to the fact that, when $\muhat$ is indivisible, $t^{d_\mu}A_\mu(t^2)$ is \cite{crawley-boevey-etal} the compactly supported Poincar\'e polynomial of $\calQ_\mu$ and thus the middle Betti number of $\calQ_\mu$ is $A_\mu(0)$.

Also, Theorems \ref{multi} and \ref{abs} imply that $\langle
\Lambda\otimes R_{\muhat},1\rangle=A_{\muhat}(q)$. This gives an
unexpected connection between the representation theory of
$\GL_n(\F_q)$ and that of comet-shaped, typically wild, quivers.
\section{Connectedness of character varieties}

The quiver variety $\calQ_{\muhat}$ is known to be connected
\cite{crawley-moment}. Here we use Theorem \ref{epolychar} to prove
the following theorem~\cite{hausel-letellier-villegas2}:

\begin{theorem} The character variety $\M_{\muhat}$ is
connected. \label{connected}\end{theorem}

Since the character variety $\M_{\muhat}$ is non-singular, the mixed
Hodge numbers $h^{i,j;k}(\M_{\muhat})$ equal zero if
$(i,j,k)\notin\{(i,j,k)|\, i\leq k,\,j\leq k,\, k\leq i+j\}$, see
\cite{Del1}. The number of connected components of $\M_{\muhat}$ is
equal to $h^{0,0;0}(\M_{\muhat})$ and $h^{0,0;k}(\M_{\muhat})=0$ if
$k>0$. Hence by Corollary \ref{curious}, we see that the number of
connected components of $\M_{\muhat}$ equals the constant term of
the $E$-polynomial $E(\M_{\muhat};q)$. To prove the theorem, we are
thus reduced to prove that the coefficient of the lowest power
$q^{\frac{1}{2}d_{\muhat}}$ of $q$ in
$\H_{\muhat}(\sqrt{q};1/\sqrt{q})$ is $1$.

We use the following expansion \cite[Lemma
5.1.5]{hausel-letellier-villegas}:

$$\sum_{\muhat\in\P^k}\frac{q\H_{\muhat}(\sqrt{q},1/\sqrt{q})}{(q-1)^2}\,m_{\muhat}=
\Log\left(\sum_{\lambda\in\P}\mathcal{H}_{\lambda}(\sqrt{q},1/\sqrt{q})\big(q^{-n(\lambda)}H_{\lambda}(q)\big)^k\prod_{i=1}^ks_{\lambda}(\x_i\y)\right)$$
where $\y=\{1,1,q^2,...\}$, $H_{\lambda}(q)$ is the hook polynomial
and $s_{\lambda}$ is the Schur symmetric function. The key-point in
 the proof of Theorem \ref{connected} for $g>0$ is the following
result \cite{hausel-letellier-villegas2}:

\begin{theorem} Given a partition $\lambda\in\P$, let
$v(\lambda)$ be the lowest power of $q$ in
$$\mathcal{A}_{\lambda}(q):=\mathcal{H}_{\lambda}(\sqrt{q},1/\sqrt{q})\big(q^{-n(\lambda)}H_{\lambda}(q)\big)^k\prod_{i=1}^k\left\langle
h_{\mu^i}(\x_i),s_{\lambda}(\x_i\y)\right\rangle.$$ If $g>0$, then the minimum of
the $v(\lambda)$ where $\lambda$ runs over the partitions of a given
size $n$, occurs only at $\lambda=(1^n)$. Moreover
$v(\lambda)=-\frac{1}{2}d_{\muhat}+1$ and the coefficient of
$q^{-\frac{1}{2}d_{\muhat}+1}$ in $\mathcal{A}_{(1^n)}(q)$ is $1$.
\label{optim}\end{theorem}

When $g=0$, Theorem \ref{optim} is known to fail is some cases. Instead we proceed with a proof which combines the use of Weyl symmetry or Katz convolution  at the middle vertex and an analogue of Theorem \ref{optim}. Here the partition $\lambda=(1^n)$ may be not the only one for which $v(\lambda)$ is minimal.  However, we show that an appropriate cancellation occurs after taking the Log.  \vspace{.2cm}

\section{Relation with Hilbert schemes on $\C^*\times\C^*$}

Put $X:=\C^*\times\C^*$ and denote by $X^{[n]}$ the Hilbert scheme
of $n$ points on $X$. We have \cite{hausel-letellier-villegas2}:

\begin{theorem}
 Assume that $g=1$ and $\muhat$ is the single partition
$\mu=(n-1,1)$. Then $X^{[n]}$ and $\M_{\muhat}$ have the same mixed
Hodge polynomial. \label
{hilbert-hodge}
\end{theorem}
The compactly supported mixed Hodge polynomial of $X^{[n]}$ is given
by the following generating function~\cite{Gottsche-Soergel}:
\begin{equation}
\label{goettsche}
1+\sum_{n\geq 1}H_c\big(X^{[n]};q,t)T^n=\prod_{n\geq
  1}\frac{(1+t^{2n+1}q^nT^n)^2}{(1-q^{n-1}t^{2n}T^n)(1-t^{2n+2}q^{n+1}
  T^n)}.
\end{equation}
The identity~\eqref{goettsche} combined with the case $g=1$ and
$\mu=(n-1,1)$ of our Conjecture~\ref{conj} becomes the following
purely combinatorial conjectural identity:
\begin{conjecture}
\begin{equation}
1+(z^2 - 1) (1-w^2)
\frac
{ \sum_\lambda \mathcal{H}_\lambda (z,w) \phi_\lambda  (z^2,w^2)T^{|\lambda|}}
{\sum_\lambda \mathcal{H}_\lambda (z,w) T^{|\lambda|}}=
\prod_{n\geq
1}\frac{(1-zwT^n)^2}{(1-z^2T^n)(1-w^2T^n)},
\label{hilbert}
\end{equation}
where $\phi(0):=0$ and if $\lambda$ is a non-zero partition
$$
\phi_\lambda (z,w) :=
 \sum_{(i,j)\in \lambda} z^{j-1} w^{i-1},
$$
where the sum runs over the boxes of $\lambda$.
\label{conj2}
\end{conjecture}

\begin{theorem} Equation \eqref{hilbert} is true in the specialization
$(z,w)\mapsto (1/\sqrt{q},\sqrt{q})$.
\end{theorem}
This theorem is a consequence of \eqref{goettsche}, Theorems
\ref{epolychar} and \ref{hilbert-hodge}; in
\cite{hausel-letellier-villegas2} we give an alternative purely
combinatorial proof. Putting $q=e^u$ yields the following
\begin{corollary}
$$
1+\sum_{n\geq 1}
\H_{\muhat}\big(e^{u/2},e^{-u/2})T^n=
\frac{1}{u}\big(e^{u/2}-e^{-u/2}\big)
\exp\left(2\sum_{k\geq2}G_k\frac{u^k}{k!}\right)
$$
where $G_k,\, k\geq 2$ are the standard Eisenstein series. In
particular, the coefficient of any power of $u$ of the left hand side
is in the ring of \emph{quasi-modular forms}, generated by the $G_k$,
$k\geq 2$ over $\Q$.
\end{corollary}
The fact that modular forms might be involved in this situation was
pointed out in \cite{vafa}, see also \cite{Goettsche} and
\cite{Bloch-Okounkov}.


\begin{thebibliography}{00}



\bibitem{Bloch-Okounkov}{\sc S. Bloch {\rm and} A. Okounkov}
The character of the infinite wedge representation.
Adv. Math. 149 (2000), no. 1, 1--60


\bibitem{crawley-moment} {\sc W. Crawley-Boevey}: Geometry of the
moment map for representations of quivers. {\em Comp. Math.} {\bf
126} (2001), 257--293.

\bibitem{crawley-boevey-etal} {\sc Crawley-Boevey, W. {\rm and} Van den Bergh, M.}: \newblock
Absolutely indecomposable representations and Kac-Moody Lie
algebras. With an appendix by Hiraku Nakajima. \newblock {\em
Invent. Math.} {\bf 155} (2004), no. 3, 537--559.


\bibitem{Del1}{\sc P. Deligne}:\newblock Th\'eorie de Hodge II. {\em
Inst. hautes Etudes Sci. Publ. Math.} {\bf  40} (1971), 5--47.


\bibitem{garsia-haiman} {\sc Garsia, A.M. {\rm and} Haiman, M.}:
  \newblock  A remarkable q,t-Catalan sequence and q-Lagrange
  inversion,
\newblock
{\em J. Algebraic Combin.} {\bf 5} (1996) no. 3, 191--244.

\bibitem{Goettsche} {\sc L. G\"ottshce}: \newblock Theta functions and Hodge
  numbers of moduli spaces of sheaves on rational surfaces.
\newblock
{\em Comm. Math. Phys.} {\bf 206},  (1999), 105--136.

\bibitem{Gottsche-Soergel} {\sc L. G\"ottsche {\rm and } W.
Soergel}:\newblock Perverse sheaves and the cohomology of Hilbert
schemes of smooth algebraic surfaces. {\em Math. Ann. }{\bf 296}
(1993), 235--245.

\bibitem{hausel2} {\sc T. Hausel}: Arithmetic harmonic analysis, Macdonald polynomials and
the topology of the Riemann-Hilbert monodromy map (with an Appendix
by E. Letellier) (in preparation)


\bibitem{hausel-letellier-villegas}
{\sc T. Hausel E. Letellier {\rm and} F. Rodriguez-Villegas},
{Arithmetic harmonic analysis on character and quiver varieties}.
arXiv:0810.2076.

\bibitem{hausel-letellier-villegas2}
{\sc T. Hausel E. Letellier {\rm and} F. Rodriguez-Villegas},
{Arithmetic harmonic analysis on character and quiver varieties II}.
In preparation.

\bibitem{hua} {\sc J.~Hua}: Counting representations of quivers over finite fields. {\em J. Algebra} {\bf 226}, (2000) 1011--1033


\bibitem{kac}{\sc V. Kac}: Root systems, representations of quivers
and invariant theory. {\em Lecture Notes in Mathematics}, {\bf vol.
996}, Springer-Verlag (1982), 74--108.


\bibitem{vafa}{\sc C. Vafa {\rm and } E. Witten}: A strong coupling test
of S-duality. {\em Nuclear Phys. B}  {\bf 431}  (1994),  no. 1-2,
3--77.




\end{thebibliography}
\end{document}